\documentclass[11pt,a4paper]{article}
\usepackage[T1]{fontenc}
\usepackage{lmodern,microtype}
\usepackage[margin=28mm]{geometry}
\usepackage{amsmath,amssymb,amsthm,mathtools}
\usepackage{booktabs,array}
\usepackage{xcolor}
\usepackage[colorlinks=true,linkcolor=blue!50!black,citecolor=blue!50!black,urlcolor=blue!50!black]{hyperref}
\hypersetup{pdftitle={Higher order MDS codes from normal rational curves},pdfauthor={}}
\newtheorem{theorem}{Theorem}[section]

\newtheorem{lemma}[theorem]{Lemma}

\theoremstyle{definition}
\newtheorem{definition}[theorem]{Definition}
\theoremstyle{remark}

\newcommand{\F}{\mathbb F}
\newcommand{\PG}{\operatorname{PG}}
\newcommand{\MDS}{\operatorname{MDS}}
\newcommand{\spn}[1]{\langle #1\rangle}
\newcommand{\cC}{\mathcal C}
\newcommand{\cD}{\mathcal D}
\newcommand{\cQ}{\mathcal Q}
\newcommand{\cH}{\mathcal H}
\newcommand{\cV}{\mathcal V}
\newcommand{\cA}{\mathcal A}
\newcommand{\ii}{\infty}
\allowdisplaybreaks[1]
\begin{document}
\title{Higher order MDS codes from normal rational curves}
\author{Francesco Pavese$^{1}$ and Ago-Erik Riet$^{2}$ \\ [0.5em]
\small $^{1}$Dipartimento di Meccanica, Matematica e Management\\
\small Politecnico di Bari, Bari, Italy\\
\small francesco.pavese@poliba.it \\ [0.3em] 
\small $^{2}$Institute of Mathematics and Statistics\\
\small University of Tartu, Tartu, Estonia \\
\small ago-erik.riet@ut.ee} 
\date{}
\maketitle
\vspace{-1.3em}
\begin{abstract}
We give explicit constructions of size at least $q$ of arcs in PG$(2,q^3)$ with no three secants intersecting in a common external point, equivalently point sets in general position on level 2, equivalently the first $[\ge q,3]$-MDS$(3)$ codes over $\F_{q^3}$. More precisely, we give geometric criteria for projecting a selected set of points defined over $\F_q$ of a twisted cubic or a normal quartic to an arc $\cA$ in $\PG(2,q^3)$ such that no three distinct secant lines to $\cA$ have a point in common outside $\cA$. Applications give $q+1$ points when $q$ is even or an odd power of $3$, and $q$ points for every odd $q$. 
\end{abstract}

\section{Introduction}

The geometry of higher-order MDS codes strengthens the usual arc condition. In dimension three, the columns of a generator matrix must not only avoid collinear triples: the lines joining three disjoint pairs of columns must also avoid a common point. This interpretation is part of the theory introduced by Brakensiek, Gopi and Makam~\cite{BGM}; see also Roth~\cite{Roth} for the related list-decoding formulation. Nagy and Weiner~\cite{NW}, apparently unaware of the work~\cite{BGM}, define the equivalent notion of general position on level 2, used to study secant size distributions.

Throughout, subspaces defined over $\F_q$ are extended to $\F_{q^3}$ by extending scalars; this extension is indicated by a subscript $q^3$. All curve points selected for projection are defined over $\F_q$. We use the algebraic normal curves
\begin{align}
 \cD(\F_q)&=\{D(t)=(1:t:t^2:t^3):t\in \F_q\}\cup\{D(\ii)=(0:0:0:1)\},\label{eq:cubic}\\
 \cC(\F_q)&=\{R(t)=(1:t:t^2:t^3:t^4):t\in \F_q\}\cup\{R(\ii)=(0:0:0:0:1)\}.\label{eq:quartic}
\end{align}

\begin{definition}\label{def:arc}
An \emph{arc} in $\PG(2,q^3)$ is a set $\cA$ of at least three points, no three of which are collinear. A \emph{secant line} to $\cA$ is a line $\spn{A,B}$ joining two distinct points $A,B\in\cA$. We consider arcs $\cA$ with the additional property
\begin{equation}\label{eq:nonconcurrence}
 \ell_1\cap\ell_2\cap\ell_3\subseteq\cA
 \quad\text{for any three distinct secant lines $\ell_1,\ell_2,\ell_3$ to $\cA$.}
\end{equation}
Thus no three distinct secant lines to $\cA$ have a point in common outside $\cA$.
\end{definition}

An $n$-point arc satisfying~\eqref{eq:nonconcurrence} gives an $[n,3]$-$\MDS(3)$ code by taking homogeneous coordinates as columns, and conversely~\cite{BGM,BDG}.

For the following comparison, let $q$ denote the order of the ambient field. To the best of our knowledge, the largest previously published explicit families of arcs satisfying~\eqref{eq:nonconcurrence}, in terms of asymptotic size, are in characteristic $7$. The construction of Brakensiek, Dhar and Gopi~\cite[Theorem~1.13 and Section~5.2]{BDG} gives at least $q^{1/3}/7$ points for the cubic extension fields covered by that construction.

In this note we obtain arcs satisfying~\eqref{eq:nonconcurrence} of size at least $q^{1/3}$ for every cube prime power $q$, in every characteristic, i.e.~gen$_2(q)\ge q^{1/3}$ in the notation of Nagy and Weiner~\cite{NW}. The size is $q^{1/3}$ in odd characteristic and $q^{1/3}+1$ in characteristic two; the latter size is also obtained when the base field has order an odd power of $3$. In the remainder, $q$ again denotes the order of the base field.

\section{Geometric criteria for selected subsets}\label{sec:criteria}

Let $r\in\{3,4\}$ and let $S$ be a subset of the points defined over $\F_q$ of the corresponding normal curve, with $|S|\ge3$. Set
\begin{equation}\label{eq:chordunion}
 \cV_{q^3}(S)=\bigcup_{\{P,Q\}\subset S,\ P\ne Q}\spn{P,Q}_{q^3}\subset\PG(r,q^3).
\end{equation}
The lines $\spn{P,Q}_{q^3}$ occurring in this union are called \emph{chords of $\cV_{q^3}(S)$}; thus each joins two distinct points of $S$, and tangents are not included. Two chords determined by four distinct points of $S$ are skew, since any four distinct points of either normal curve are linearly independent. Consequently, two distinct chords can meet only at a common point of $S$. The set $\cV_{q^3}(S)$ is a finite union of full $\F_{q^3}$-lines. It is not the algebraic secant variety of the curve: chords involving points outside $S$ and tangent lines play no role in~\eqref{eq:chordunion}.

Let $\Lambda$ be a projective $(r-3)$-space over $\F_{q^3}$. The quotient by its vector space defines a projection to $\PG(2,q^3)$.

\begin{theorem}\label{thm:intersection}
Projection from $\Lambda$ is defined and injective on $S$ and has image an arc satisfying~\eqref{eq:nonconcurrence} if and only if
\begin{equation}\label{eq:intersectioncriterion}
 |\mathcal F\cap\cV_{q^3}(S)|\le2
 \quad\text{for every $(r-2)$-space $\mathcal F\supset\Lambda$.}
\end{equation}
\end{theorem}

\begin{proof}
Let $\varphi$ denote projection from $\Lambda$. Assume~\eqref{eq:intersectioncriterion}. If $\Lambda$ met a chord $\ell$ of $\cV_{q^3}(S)$, the span $\spn{\Lambda,\ell}$ would be contained in an $(r-2)$-space through $\Lambda$. This $(r-2)$-space would contain the whole chord, contrary to~\eqref{eq:intersectioncriterion}. Consequently the projection is defined and injective on $S$, and maps every chord of $\cV_{q^3}(S)$ isomorphically onto a line.

Take three distinct points $A,B,C\in S$ and let $\pi=\spn{A,B,C}_{q^3}$, a plane. Their images are collinear exactly when $\pi\cap\Lambda\ne\varnothing$. Since $\Lambda$ meets no chord of $\cV_{q^3}(S)$, $\pi\cap\Lambda$ is a single point $N$, not on a side of the triangle $ABC$. Choose a line $m\subset \pi$ through $N$ avoiding the three vertices. It meets the three sides in distinct points. The space $\spn{\Lambda,m}$ is a line through the centre or plane through the centre and contains these three points of $\cV_{q^3}(S)$, a contradiction. The choice of $m$ is possible since $|\F_{q^3}|\ge8$.

If the three projected secants of three disjoint chords of $\cV_{q^3}(S)$ were concurrent, their three intersection points with the corresponding $(r-2)$-space through $\Lambda$ would be distinct, again contradicting~\eqref{eq:intersectioncriterion}. These intersections are distinct because the three chords are pairwise skew, by the preceding observation.

Conversely, suppose the projection is defined and injective on $S$, its image is an arc satisfying~\eqref{eq:nonconcurrence}, and an $(r-2)$-space $\mathcal F\supset\Lambda$ contains distinct points $U_1,U_2,U_3\in\cV_{q^3}(S)$. Choose chords $\ell_i$ of $\cV_{q^3}(S)$ containing $U_i$. None is contained in $\mathcal F$, since the two points of $S$ defining it have different images. Hence each meets $\mathcal F$ in exactly one point, and the $\ell_i$ are distinct. Their images are three distinct lines through the same image point $P$; equality of two of these lines would violate the arc condition. If two chords shared a point $A\in S$, their projected lines would meet at the image of $A$, so $P=\varphi(A)$. Injectivity on each chord would then force the two corresponding $U_i$ to equal $A$, a contradiction. The three chords therefore have six distinct points of $S$ and give a forbidden concurrence. This proves the converse.
\end{proof}

\subsection{Twisted cubics: trisecant planes and hyperbolic quadrics}\label{sec:cubicapps}

Fix $S\subseteq\cD(\F_q)$. Let $\cH(S)$ be the set of planes spanned by three points of $S$, after extension over $\F_{q^3}$. Define $\cQ(S)$ to be the set of nonsingular hyperbolic quadrics defined over $\F_q$ that contain three pairwise skew chords of $\cV_{q^3}(S)$ in one ruling, after extension to $\F_{q^3}$. 

\begin{theorem}\label{thm:cubiccriterion}
Let $S\subseteq\cD(\F_q)$ and $O\in\PG(3,q^3)$. Projection from $O$ is defined and injective on $S$ and has image an arc satisfying~\eqref{eq:nonconcurrence} if and only if
\begin{equation}\label{eq:badcubic}
 O\notin
 \bigcup_{H\in\cH(S)}H\ \cup\
 \bigcup_{Q\in\cQ(S)}Q.
\end{equation}
\end{theorem}
\begin{proof}
Avoiding $\cH(S)$ is equivalent to being defined and injective on $S$ and having no collinear triple in the image. Indeed, a selected pair is contained in a selected triple, and three projected points are collinear exactly when their spanning plane contains $O$.

Assume this condition. The three chords of six distinct curve points are pairwise skew. By the standard regulus property of three pairwise skew lines, a line through $O$ meets these three chords if and only if $O$ lies on their hyperbolic quadric. If $O$ is on the quadric, its opposite-ruling line provides such a line, with three distinct intersections. The point $O$ cannot lie on one of the chords of $\cV_{q^3}(S)$, by the first condition. This proves~\eqref{eq:badcubic}.
\end{proof}

Let $\sigma$ be the Frobenius collineation of $\PG(3,q^3)$ induced by the map
\[
 (x_0,\ldots,x_3)\longmapsto(x_0^q,\ldots,x_3^q).
\]
Every quadric defined over $\F_q$ is stabilized by $\sigma$; the same holds for every plane defined over $\F_q$.

\subsubsection{Characteristic two}

\begin{theorem}\label{thm:even}
Let $q$ be even and $\lambda\in \F_{q^3}\setminus \F_q$. The point
\[
 O_\lambda=(0:1:0:\lambda)
\]
defines a projection that is defined and injective on all $q+1$ points of $\cD(\F_q)$, with image an arc satisfying~\eqref{eq:nonconcurrence}. The quotient
\[
 (X_0:X_1:X_2:X_3)\longmapsto
 (X_3+\lambda X_1:X_0:X_2+\lambda X_0)
\]
has image
\begin{equation}\label{eq:evenpoints}
 \cA_\lambda=\{(t^3+\lambda t:1:t^2+\lambda):t\in \F_q\}
                  \cup\{(1:0:0)\}.
\end{equation}
\end{theorem}
\begin{proof}
Put $S=\cD(\F_q)$. The point $O_\lambda$ lies in $\ell\setminus\PG(3,q)$, where $\ell:X_0=X_2=0$. Any plane or quadric defined over $\F_q$ containing $O_\lambda$ is stabilized by $\sigma$, so it contains its three distinct conjugates on $\ell$, and hence contains $\ell$. A plane through $\ell$ contains at most two points of $\cD(\F_q)$. Thus no plane spanned by three points of $S$ contains $O_\lambda$.

Suppose a quadric $Q$ of $\cQ(S)$ contains $O_\lambda$. Thus $Q$ contains $\ell$. The line $\ell$ meets no chord with two distinct points of $S\setminus\{D(\ii)\}$. If one of the three defining chords of $Q$ contains $D(\ii)$, the line $\ell$ meets exactly one of the three ruling lines. This is impossible for a line on a nonsingular hyperbolic quadric: it meets either all three or none, unless it is one of them; here $\ell$ is not a chord of $\cV_{q^3}(S)$.

Otherwise $Q$ contains six curve points and $D(\ii)\in\ell$. By Hirschfeld~\cite[Corollary~2, p.~237]{Hirschfeld}, a hyperbolic quadric containing seven distinct points of $\cD$ contains the entire curve. Hence $Q$ contains $\cD$. The quadratic ideal of the twisted cubic is generated by
\begin{equation}\label{eq:ideal}
 X_0X_2-X_1^2,\quad X_0X_3-X_1X_2,\quad X_1X_3-X_2^2.
\end{equation}
Containment of $\ell$ leaves only $X_0X_3+X_1X_2=0$. Its polar form on $D(u),D(v)$ equals $(u+v)^3\ne0$, so it contains no chord of $\cV_{q^3}(S)$ joining two points of $S\setminus\{D(\ii)\}$, a contradiction. Theorem~\ref{thm:cubiccriterion} applies. The quotient formula verifies~\eqref{eq:evenpoints} directly.
\end{proof}

\subsubsection{Odd powers of three}

\begin{theorem}\label{thm:three}
Let $q=3^{2m+1}$, $m\ge0$, and $\gamma\in \F_{q^3}\setminus \F_q$. The point
\[
 O_\gamma=(0:1:-\gamma:0)
\]
defines a projection that is defined and injective on all of $\cD(\F_q)$, with image an arc satisfying~\eqref{eq:nonconcurrence}. The quotient $(X_0:X_3:X_2+\gamma X_1)$ gives
\begin{equation}\label{eq:threepoints}
 \cA_\gamma=\{(1:t^3:t^2+\gamma t):t\in \F_q\}
                      \cup\{(0:1:0)\}.
\end{equation}
\end{theorem}
\begin{proof}
Put $S=\cD(\F_q)$ and $\ell:X_0=X_3=0$.
First, if a quadric $Q$ of $\cQ(S)$ contains $O_\gamma$, then $Q$ is stabilized by $\sigma$, since it is defined over $\F_q$. It therefore contains the three distinct points $O_\gamma$, $O_{\gamma^q}$ and $O_{\gamma^{q^2}}$ of $\ell$, and hence contains $\ell$.

In characteristic three the osculating planes are $X_3-t^3X_0=0$ and $X_0=0$. Their common axis is $\ell$, preserved by the natural stabilizer $\operatorname{PGL}(2,\F_q)$. A plane through $\ell$ contains at most one cubic point defined over $\F_q$. A plane defined over $\F_q$ containing $O_\gamma$ is stabilized by $\sigma$, and thus contains $O_{\gamma^q}$ as well. These two points span $\ell$. Therefore no plane spanned by three points of $S$ contains $O_\gamma$.

Suppose a quadric $Q$ of $\cQ(S)$ contains $\ell$. Using the stabilizer, normalize one of its three defining chords to $r=\spn{D(0),D(\ii)}$; the axis $\ell$ is unchanged. Let $\spn{D(u),D(v)}$ be a further defining chord, where $u,v\in\F_q^*$ are distinct, and put $s=u+v$, $p=uv$. A quadric containing $\ell$ and $r$ has equation
\[
 A X_0X_1+B X_0X_2+C X_1X_3+D X_2X_3=0.
\]
Substitution of $\alpha D(u)+\beta D(v)$ and comparison of the coefficients of $\alpha^2$, $\beta^2$ and $\alpha\beta$ give
\begin{align*}
 A+Bu+Cu^3+Du^4&=0,\\
 A+Bv+Cv^3+Dv^4&=0,\\
 As+B(s^2+p)+Cp(s^2+p)+Dp^2s&=0.
\end{align*}
Solving these three linear equations gives
\[
 (A,B,C,D)=\eta(-p^2,ps,-s,1),\qquad \eta\ne0.
\]
Thus the quadric containing $\ell$, $r$ and $\spn{D(u),D(v)}$ is
\begin{equation}\label{eq:threequadric}
 -p^2X_0X_1+psX_0X_2-sX_1X_3+X_2X_3=0.
\end{equation}
It is nonsingular, since the matrix of the terms joining the coordinate pairs $(X_0,X_3)$ and $(X_1,X_2)$ has determinant
\[
 \det\begin{pmatrix}-p^2&ps\\-s&1\end{pmatrix}
 =p(s^2-p)=p(u-v)^2\ne0.
\]

Write the third defining chord as $\spn{D(w),D(z)}$, and put $s'=w+z$, $p'=wz$. Applying the same calculation to this chord and comparing~\eqref{eq:threequadric}, whose $X_2X_3$ coefficient is normalized to $1$, gives $s'=s$, $(p')^2=p^2$ and $p's=ps$. If $p'=p$, the two pairs coincide. Otherwise $p'=-p$ and $s=0$, so the pairs are $\{u,-u\}$ and $\{w,-w\}$ with $w^2=-u^2$. This contradicts the fact that $-1$ is a nonsquare in $\F_q$ when $q=3^{2m+1}$.

Hence no quadric of $\cQ(S)$ contains $\ell$, and the initial observation shows that none contains $O_\gamma$. Theorem~\ref{thm:cubiccriterion} completes the proof.
\end{proof}

\subsection{Quartics: projection from a line \texorpdfstring{$\Lambda$}{Lambda}}

For $S\subseteq\cC(\F_q)$, set
\begin{equation}\label{eq:badquartic}
 \begin{aligned}
 \mathcal B(S)=\{\Lambda:\;&\Lambda\text{ is a line in }\PG(4,q^3),\\
 &\exists\text{ a plane }\Pi\supset\Lambda
       \text{ with }|\Pi\cap\cV_{q^3}(S)|\ge3\}.
 \end{aligned}
\end{equation}
By Theorem~\ref{thm:intersection}, projection from $\Lambda$ is defined and injective on $S$ and has image an arc satisfying~\eqref{eq:nonconcurrence} precisely when $\Lambda\notin\mathcal B(S)$. Equivalently, $\Lambda$ must avoid every plane spanned by three points of $S$ and admit no plane through $\Lambda$ meeting three chords determined by six distinct points of $S$. The first requirement already excludes intersections of $\Lambda$ with chords of $\cV_{q^3}(S)$.

\subsubsection{A construction for odd \texorpdfstring{$q$}{q}}\label{sec:odd}

\begin{lemma}\label{lem:existence-a}
Let $q$ be odd and let $\delta\in\F_q$ be a nonsquare. There exists $a\in\F_q$ such that $a^2-\delta$ is a nonsquare. There are exactly $(q+1)/2$ such elements $a$. If $q\equiv1\pmod4$, one may take $a=0$.
\end{lemma}
\begin{proof}
Let $\chi$ be the quadratic character of $\F_q$. The equation $y^2=a^2-\delta$ has $q-1$ solutions $(a,y)$, since $(a-y)(a+y)=\delta$. Thus $\sum_{a\in \F_q}\chi(a^2-\delta)=-1$. There are no zero terms, so exactly $(q+1)/2$ of them equal $-1$. Finally, $-\delta$ is a nonsquare when $q\equiv1\pmod4$, giving $a=0$.
\end{proof}

Throughout this subsection, let $q$ be odd, choose $a,\delta\in\F_q$ such that both $\delta$ and $a^2-\delta$ are nonsquares, and let $\gamma\in\F_{q^3}\setminus\F_q$. Put
\[
 S=\{R(t):t\in\F_q\}\subset\cC(\F_q).
\]
In particular, $R(\ii)$ is not selected. Define
\begin{equation}\label{eq:centreodd}
 \Lambda_{a,\delta,\gamma}=
 \spn{(1,0,\delta,0,\delta^2-a\gamma),\
       (0,1,0,\delta,-\gamma)}_{q^3}.
\end{equation}
For brevity, write $\Lambda=\Lambda_{a,\delta,\gamma}$. The plane
\begin{equation}\label{eq:auxiliaryplane}
 E:\quad X_2=\delta X_0,\qquad X_3=\delta X_1
\end{equation}
is defined over $\F_q$ and contains $\Lambda$. Moreover,
\begin{equation}\label{eq:rationalcentrepoint}
 \Lambda\cap\PG(4,q)=\{A\},\qquad
 A=(1:-a:\delta:-a\delta:\delta^2).
\end{equation}
Indeed, a point of $\Lambda$ has coordinates
\[
 (x:y:\delta x:\delta y:\delta^2x-\gamma(ax+y)).
\]
If this point is defined over $\F_q$, its coordinates may be scaled so that $x,y\in\F_q$. Since $\gamma\notin\F_q$, its last coordinate then forces $ax+y=0$, which gives $A$.

We shall use projection from $A$ as an intermediate step. The map
\begin{equation}\label{eq:projectionA}
 \rho:(X_0:\cdots:X_4)\longmapsto(X_1+aX_0:X_2-\delta X_0:X_3-\delta X_1:X_4-\delta X_2)
\end{equation}
has vector rank four and projective kernel $A$. In its target, with coordinates $(Y_0:Y_1:Y_2:Y_3)$, we have
\begin{equation}\label{eq:auxiliarylinepoint}
 e=\rho(E\setminus\{A\}):Y_1=Y_2=0,
 \qquad
 O=\rho(\Lambda\setminus\{A\})=(1:0:0:-\gamma).
\end{equation}
We use $\sigma$ for the Frobenius collineation on this copy of $\PG(3,q^3)$ as well. The three points $O,O^\sigma,O^{\sigma^2}$ are distinct and lie on the line $e$, which is defined over $\F_q$.

For later use, $A$ lies in no plane spanned by three distinct points of $S$. To see this, work temporarily over $\F_{q^6}$ and choose $v\in\F_{q^6}$ with $v^2=\delta$. Since $\delta$ is a nonsquare in $\F_q$, we have $v\in\F_{q^2}\setminus\F_q$ and $v\notin\F_{q^3}$. Then
\begin{equation}\label{eq:Aimaginarychord}
 A=\frac{v-a}{2v}R(v)+\frac{v+a}{2v}R(-v),
\end{equation}
where the equality uses the displayed vector representatives. Both coefficients are nonzero. Any three distinct points of $S$, together with $R(v)$ and $R(-v)$, are five distinct points of the normal quartic and are linearly independent by the Vandermonde determinant. Equation~\eqref{eq:Aimaginarychord} therefore excludes $A$ from the plane spanned by the three selected points. It follows that $\rho$ is defined and injective on $S$, its image has no three collinear points, and it maps every chord of $\cV_{q^3}(S)$ isomorphically onto a line.

\begin{lemma}\label{lem:quartic-transversal}
No line through $A$ meets three chords determined by six distinct points of $S$.
\end{lemma}
\begin{proof}
Suppose that a line $r$ through $A$ meets
\[
 r_i=\spn{R(t_i),R(u_i)}_{q^3},\qquad 1\le i\le3,
\]
where the six parameters are distinct. The chords $r_i$ are pairwise skew. Thus their intersections $Z_i=r\cap r_i$ are distinct. None of the $Z_i$ is a point of $S$ defining its chord: otherwise the collinearity of the $Z_i$ would give a nonzero linear relation among at most five distinct points of the quartic, which is impossible.

The $5\times6$ matrix whose columns are the six curve vectors has rank five. Its kernel is one-dimensional, is defined over $\F_q$, and has no zero coordinate, since every five columns are independent. Choose its relation in the form
\begin{equation}\label{eq:sixpointrelation}
 \sum_{i=1}^3\bigl(c_iR(t_i)+d_iR(u_i)\bigr)=0,
 \qquad c_i,d_i\in\F_q^*.
\end{equation}
Collinearity of the $Z_i$ also gives a relation among these six curve vectors. Uniqueness of~\eqref{eq:sixpointrelation} implies
\[
 Z_i=\bigl[c_iR(t_i)+d_iR(u_i)\bigr],\qquad 1\le i\le3,
\]
where brackets denote the projective point represented by a nonzero vector. Consequently each $Z_i$, and hence $r$, is defined over $\F_q$.

For coordinates $(X_0,\ldots,X_4)$, define
\begin{equation}\label{eq:hankelquartic}
 H(X_0,\ldots,X_4)=
 \begin{pmatrix}
 X_0&X_1&X_2\\
 X_1&X_2&X_3\\
 X_2&X_3&X_4
 \end{pmatrix}.
\end{equation}
If $\mathbf v(t)=(1,t,t^2)^{\mathsf T}$, then
\[
 H(R(t))=\mathbf v(t)\mathbf v(t)^{\mathsf T}.
\]
Thus the homogeneous cubic $\det H(X_0,\ldots,X_4)$ vanishes on every chord. It also vanishes at $A$, since
\[
 H(A)=
 \begin{pmatrix}
 1&-a&\delta\\
 -a&\delta&-a\delta\\
 \delta&-a\delta&\delta^2
 \end{pmatrix},
 \qquad
 \ker H(A)=\spn{(-\delta,0,1)^{\mathsf T}}.
\]
The point $A$ is not on any chord, by the observation preceding the lemma. Hence the restriction of $\det H$ to $r$ has the four distinct zeros $A,Z_1,Z_2,Z_3$. A nonzero homogeneous cubic on a projective line has at most three distinct zeros, so this restriction is identically zero.

Take the basis $(1,0,0)^{\mathsf T},(0,1,0)^{\mathsf T},(-\delta,0,1)^{\mathsf T}$ of $\F_q^3$. In this basis, the symmetric matrix $H(A)$ is congruent to $\operatorname{diag}(B,0)$, where
\[
 B=\begin{pmatrix}1&-a\\-a&\delta\end{pmatrix}.
\]
The binary quadratic form associated with $B$ is anisotropic over $\F_q$: indeed,
\[
 x^2-2axy+\delta y^2=(x-ay)^2-(a^2-\delta)y^2
\]
vanishes only when $x=y=0$, because $a^2-\delta$ is a nonsquare.

Choose a rational vector $T$ such that $A,T$ span $r$, and write the matrix congruent to $H(T)$ as
\[
 \begin{pmatrix}C&b\\b^{\mathsf T}&c\end{pmatrix},
 \qquad C\in\operatorname{Mat}_{2\times2}(\F_q),\quad
 b\in\F_q^2,\quad c\in\F_q.
\]
The polynomial identity just obtained is
\begin{equation}\label{eq:singularpencil}
 \det\begin{pmatrix}B+zC&zb\\zb^{\mathsf T}&zc\end{pmatrix}=0
 \quad\text{in }\F_q[z].
\end{equation}
Its coefficient of $z$ is $c\det B$, so $c=0$. With $c=0$, its coefficient of $z^2$ is $-b^{\mathsf T}\operatorname{adj}(B)b$. The form associated with
\[
 \operatorname{adj}(B)=\begin{pmatrix}\delta&a\\a&1\end{pmatrix}
\]
is also anisotropic, since
\[
 \delta x^2+2axy+y^2=(y+ax)^2-(a^2-\delta)x^2.
\]
Therefore $b=0$. It follows that every matrix $H(X_0,\ldots,X_4)$ associated with a point $(X_0:\cdots:X_4)\in r$ annihilates $(-\delta,0,1)^{\mathsf T}$.

On the other hand, for the rational representative $Z_i=c_iR(t_i)+d_iR(u_i)$ we obtain
\[
 H(Z_i)\begin{pmatrix}-\delta\\0\\1\end{pmatrix}
 =c_i(t_i^2-\delta)\mathbf v(t_i)
  +d_i(u_i^2-\delta)\mathbf v(u_i).
\]
The two vectors on the right are linearly independent, and both coefficients are nonzero. This contradicts the common kernel and proves the lemma.
\end{proof}

\begin{samepage}
\begin{lemma}\label{lem:quartic-auxiliarygeometry}
In the quotient by $A$, the following statements hold.
\begin{enumerate}
\item No plane defined over $\F_q$ containing $e$ contains the image of a chord of $\cV_{q^3}(S)$.
\item No nonsingular hyperbolic quadric defined over $\F_q$ contains $e$ and the images of three chords determined by six distinct points of $S$.
\end{enumerate}
\end{lemma}
\end{samepage}
\begin{proof}
For the first assertion, a plane through $e$ lifts to a hyperplane through $E$. Its equation is $b(X_2-\delta X_0)+c(X_3-\delta X_1)=0$, for some $b,c\in\F_q$ not both zero. On the quartic this becomes
\[
 (t^2-\delta)(b+ct)=0.
\]
Since $\delta$ is a nonsquare, this equation has at most one solution in $\F_q$. The hyperplane therefore cannot contain two distinct points of $S$ defining a chord.

For the second assertion, suppose that such a quadric $\overline Q$ exists, and let $Q$ be its inverse image under~\eqref{eq:projectionA}. This inverse image is a quadric singular at $A$, containing $E$ and six distinct points of $S$. Since $Q$ contains $E$, its equation can be written
\begin{equation}\label{eq:quadricthroughE}
 \begin{aligned}
 Q(X_0,\ldots,X_4)={}&(X_2-\delta X_0)L(X_0,\ldots,X_4)\\
 &+(X_3-\delta X_1)M(X_0,\ldots,X_4),
 \end{aligned}
\end{equation}
where $L,M$ are linear forms over $\F_q$. In particular, there is no $X_4^2$ term, so the polynomial $Q(R(t))$ has degree at most seven. Its six selected roots in $\F_q$ are distinct from $v$ and $-v$, where $v^2=\delta$. Both $R(v)$ and $R(-v)$ lie in $E$, so these are two further roots. Thus $Q(R(t))$ has eight distinct roots and is identically zero.

We determine $Q$ explicitly. Substituting into~\eqref{eq:quadricthroughE} gives
\[
 0=Q(R(t))=(t^2-\delta)\bigl(L(R(t))+tM(R(t))\bigr).
\]
Consequently $L(R(t))+tM(R(t))$ is the zero polynomial. Comparing its coefficients shows that, for some $m_0,m_1,m_2,m_3\in\F_q$,
\begin{align*}
 M&=m_0X_0+m_1X_1+m_2X_2+m_3X_3,\\
 L&=-m_0X_1-m_1X_2-m_2X_3-m_3X_4.
\end{align*}
We now use that $Q$ is singular at $A$. For a homogeneous quadratic form, this means that translating the coordinates to a vector representative of $A$ produces no terms of total degree one. Thus
\[
 Q(1+y_0,-a+y_1,\delta+y_2,-a\delta+y_3,\delta^2+y_4)
\]
has no linear terms in $y_0,\ldots,y_4$. In particular, the coefficient of $z$ must vanish in each of the two restrictions obtained by varying only $X_2$ or only $X_3$.

Using~\eqref{eq:quadricthroughE} and the displayed expressions for $L$ and $M$, these restrictions are
\begin{align*}
 Q(1,-a,\delta+z,-a\delta,\delta^2)
   &=z\bigl(L(A)-m_1z\bigr),\\
 Q(1,-a,\delta,-a\delta+z,\delta^2)
   &=z\bigl(M(A)+m_3z\bigr).
\end{align*}
Indeed, on the first line $X_2-\delta X_0=z$ and $X_3-\delta X_1=0$, while on the second line these two values are $0$ and $z$, respectively. Singularity at $A$ therefore forces $L(A)=M(A)=0$. Substituting the coordinates of $A$ gives
\begin{align*}
 (m_0+\delta m_2)-a(m_1+\delta m_3)&=0,\\
 a(m_0+\delta m_2)-\delta(m_1+\delta m_3)&=0.
\end{align*}
As $a^2\ne\delta$, we obtain $m_0=-\delta m_2$ and $m_1=-\delta m_3$. Hence
\begin{align*}
 M&=m_2(X_2-\delta X_0)+m_3(X_3-\delta X_1),\\
 L&=-m_2(X_3-\delta X_1)-m_3(X_4-\delta X_2),
\end{align*}
and~\eqref{eq:quadricthroughE} reduces to
\begin{equation}\label{eq:uniquerankthreequadric}
 Q(X_0,\ldots,X_4)=m_3\bigl((X_3-\delta X_1)^2-(X_2-\delta X_0)(X_4-\delta X_2)\bigr).
\end{equation}
Since $Q$ is a nonzero quadratic form, $m_3\ne0$. In the quotient coordinates this is $Y_2^2-Y_1Y_3=0$, a quadric of rank three with singular point $(1:0:0:0)$. It cannot be the nonsingular hyperbolic quadric $\overline Q$, proving the second assertion.
\end{proof}

\begin{theorem}\label{thm:odd}
Let $q$ be odd and let $a,\delta\in\F_q$ be such that both $\delta$ and $a^2-\delta$ are nonsquares. Let $\gamma\in\F_{q^3}\setminus\F_q$.
The line $\Lambda_{a,\delta,\gamma}$ does not belong to $\mathcal B(S)$. Consequently, projection from this line is defined and injective on $S$ and has image an arc satisfying~\eqref{eq:nonconcurrence}. Its image consists of the $q$ points
\begin{equation}\label{eq:oddpoints}
 \cA_{a,\delta,\gamma}=
 \left\{\left(t:t^2+\gamma\frac{t+a}{t^2-\delta}:1\right):t\in\F_q\right\}.
\end{equation}
\end{theorem}
\begin{proof}
We first exclude intersections with planes spanned by three points of $S$, and then exclude planes through $\Lambda$ meeting three chords determined by six distinct points of $S$.

\emph{Planes spanned by three points of $S$.}
The linear map
\begin{equation}\label{eq:quarticquotient}
 (X_0,\ldots,X_4)\longmapsto
 (X_3-\delta X_1,\ X_4-\delta X_2+\gamma(aX_0+X_1),\ X_2-\delta X_0)
\end{equation}
has rank three: the columns corresponding to $X_3,X_4,X_2$ have determinant one. Its kernel is precisely the vector space defining $\Lambda$. At $R(t)$ its value is
\[
 \bigl(t(t^2-\delta),\ t^2(t^2-\delta)+\gamma(t+a),\ t^2-\delta\bigr).
\]
The last coordinate is nonzero for every $t\in\F_q$. Dividing by it gives~\eqref{eq:oddpoints}, and distinct first affine coordinates prove injectivity.

For three distinct parameters, the determinant of their three image vectors is a polynomial in $\gamma$ of degree at most one, with coefficients in $\F_q$. Its constant term is the nonzero Vandermonde determinant of the vectors $(t,t^2,1)$. Since $[\F_q(\gamma):\F_q]=3$, the determinant cannot vanish. The three image vectors are therefore independent. Equivalently, their original plane is disjoint from $\Lambda$. In particular, $\Lambda$ meets no selected chord.

\emph{Reduction to six distinct points of $S$.}
Suppose, contrary to the assertion, that a plane $\Pi\supset\Lambda$ contains three distinct points $Z_1,Z_2,Z_3$ of $\cV_{q^3}(S)$. Choose a chord $r_i$ containing each $Z_i$. No chord can be contained in $\Pi$, because two lines in a projective plane intersect, whereas each chord is disjoint from $\Lambda$. Thus each chord meets $\Pi$ in at most one point, so the $r_i$ are distinct.

The map~\eqref{eq:quarticquotient} sends $\Pi\setminus\Lambda$ to a single projective point. If two of the chords shared a point $R(t)$ of $S$, their image lines would be distinct, by the preceding noncollinearity argument, and would intersect precisely at the image of $R(t)$. Since the quotient is injective on each chord, their intersections with $\Pi$ would both equal $R(t)$, contrary to the choice of distinct $Z_i$. Hence the three chords are determined by six distinct points of $S$.

\emph{Projection of $\Pi$ from $A$.}
Let $\overline r_i=\rho(r_i)$, and let $m=\rho(\Pi\setminus\{A\})$. Then $m$ is a line through $O$ meeting all three $\overline r_i$. The lines $\overline r_i$ are defined over $\F_q$ and are distinct. Indeed, equality of two of them would put their four distinct points of $S$ in a single plane through $A$, whereas four distinct points of the normal quartic span a $3$-space.

The three lines $\overline r_i$ are not concurrent. Otherwise the inverse image of their common point would be a line through $A$ meeting all three original chords, contrary to Lemma~\ref{lem:quartic-transversal}.

Suppose that two of them, say $\overline r_1$ and $\overline r_2$, intersect at $P$. Their intersection point and their spanning plane are defined over $\F_q$. If $m$ does not pass through $P$, its intersections with these two lines are distinct, so $m$ lies in their spanning plane. If $m$ passes through $P$, then $P\notin\overline r_3$, since the three lines are not concurrent. In that case $m$ lies in the rational plane $\spn{P,\overline r_3}_{q^3}$. In either case, there is a plane defined over $\F_q$ containing $O$ and one of the lines $\overline r_i$.

Such a plane is stabilized by $\sigma$ and contains $O^\sigma$ as well as $O$. These two distinct points span $e$. The plane therefore contains $e$, contrary to the first assertion of Lemma~\ref{lem:quartic-auxiliarygeometry}. This proves that the three $\overline r_i$ must be pairwise skew.

\emph{The hyperbolic quadric determined by the three lines.}
Three pairwise skew lines in $\PG(3,q)$ determine a unique nonsingular hyperbolic quadric $\overline Q$ having them in one ruling. Its extension to $\F_{q^3}$ contains every line meeting all three, in the opposite ruling. Thus $m\subset\overline Q$, and in particular $O\in\overline Q$. Since the three defining lines are rational and the quadric is unique, $\overline Q$ is defined over $\F_q$ and is stabilized by $\sigma$.

Consequently $\overline Q$ contains the three distinct points $O,O^\sigma,O^{\sigma^2}$ of $e$. The restriction of a quadratic equation to a line cannot have three distinct zeros unless it is identically zero, so $e\subset\overline Q$. This contradicts the second assertion of Lemma~\ref{lem:quartic-auxiliarygeometry}.

No such plane $\Pi$ exists. Hence $\Lambda\notin\mathcal B(S)$, and Theorem~\ref{thm:intersection} proves the assertion. The quotient computation above gives the displayed coordinates of the image.
\end{proof}

\section{Conclusion}

Projections of twisted cubics and normal quartics give explicit arcs in $\PG(2,q^3)$ satisfying~\eqref{eq:nonconcurrence}, of size $q$ for every odd $q$ and $q+1$ when $q$ is even or an odd power of $3$. These constructions yield corresponding $[n,3]$-$\MDS(3)$ codes over $\F_{q^3}$ and a geometric description of the projection centres that must be excluded. The obtained bound $n\ge q$ in $\PG(2, q^3)$ improves the leading constant in the bound $q/7$ supplied by the construction of Brakensiek, Dhar and Gopi~\cite[Theorem~1.13 and Section~5.2]{BDG} for its applicable cubic extension fields in characteristic $7$, while extending the construction to every characteristic and every cube prime power. 

\bigskip
{\footnotesize
\noindent\textit{Acknowledgments.}
The work was supported in part by Estonian Research Council grant PRG2531 and by the Italian National Group for Algebraic and Geometric Structures and their Applications (GNSAGA--INdAM).}

\end{document}